\documentclass[10pt]{article} 

\usepackage{graphicx}
\usepackage{latexsym}
\usepackage[latin1]{inputenc} 
\usepackage{amsmath,amssymb,amsthm}
\usepackage{algorithm,algorithmic}
\usepackage{mathrsfs}
\usepackage{bm}
\usepackage{todonotes}
\usepackage{appendix} 

\usepackage{enumerate}
\usepackage{url} 
\usepackage{setspace}

\newtheorem{theorem}{Theorem}[section]
\newtheorem{lemma}[theorem]{Lemma}
\newtheorem{proposition}[theorem]{Proposition}

\newcounter{excntr}[section]

\DeclareMathOperator{\diag}{diag}

\DeclareMathOperator{\tr}{tr}

\DeclareMathOperator{\Vol}{Vol}

\newcommand{\PsiCut}{\Psi^{\mathrm{cut}}}
\newcommand{\PhiCut}{\Phi^{\mathrm{cut}}}

\newcommand{\cC}{\mathcal{C}}

\newcommand{\cG}{\mathcal{G}}

\newcommand{\cI}{\mathcal{I}}

\newcommand{\cO}{\mathcal{O}}
\newcommand{\cR}{\mathcal{R}}

\newcommand{\cOI}{\mathcal{O_\cI}}
\newcommand{\cOIplus}{\mathcal{O_\cI^+}}

\newcommand{\RR}{\mathbb{R}}

\newcommand{\x}{\times}

\begin{document}

{\title{Multiway Spectral Graph Partitioning:  Cut Functions,  Cheeger
     Inequalities,  and a Simple Algorithm}
   \author{Lars Eld\'en\\
     orcid.org/0000-0003-2281-856X\\
    Department of Mathematics, Link\"oping University, Sweden}
\maketitle
} 

\begin{abstract}
  The problem of multiway partitioning of an undirected graph is
  considered. A spectral method is used, where the $k>2$ largest
  eigenvalues of the normalized adjacency matrix (equivalently, the
  $k$ smallest eigenvalues of the normalized graph Laplacian) are computed.  It is
  shown that the information necessary for partitioning is contained
  in the subspace spanned by the $k$ eigenvectors. The partitioning is
  encoded in a matrix $\Psi$ in indicator form, which is computed by
  approximating the eigenvector matrix by a product of $\Psi$ and an
  orthogonal matrix. A measure of the distance of a graph to being
  $k$-partitionable is defined, as well as two cut (cost) functions,
  for which Cheeger inequalities are proved; thus the relation between
  the eigenvalue and partitioning problems is established. Numerical
  examples are given that demonstrate that the partitioning algorithm
  is efficient and robust.
\end{abstract}



\noindent\textbf{Running title:} Multiway Spectral Graph Partitioning

\section{Introduction}\label{sec:intro}

Let $\cG$ be an undirected graph on $n$ nodes. A key problem in many
applications (see e.g. the surveys \cite{lux07,fohr16}) is to
partition the graph into $k$ subgraphs, which are internally 
well connected, and where each subgraph is only loosely connected to
the others. Spectral 2-partitioning (partitioning into two
subgraphs, i.e. $k=2$) is based on early results by Fiedler
\cite{fied73}. This is 
a standard method for clustering that  has the advantage over many
other clustering methods in that it has a  well developed mathematical
theory, see  \cite{chung97}.
The subject of the present paper is \emph{multiway spectral
  partitioning}, or \emph{spectral $k$-partitioning},  which is the
partitioning into $k>2$ subgraphs at 
the same time. 

Traditionally, spectral partitioning is motivated by first defining a
cost function (we will call it a \emph{cut function}) for
partitioning the connected graph $\cG$ into subgraphs, 
see \cite[Section 5]{lux07}. Minimization of the cost function leads
to an optimization problem 
over all possible partitionings of the graph, which is
an NP-hard problem. Spectral partitioning is a way of solving relaxed
versions of such a problem, where the optimization is done over vectors
in $\RR^n$, and the problem becomes that of finding the $k$ smallest
eigenvalues of the  Laplacian of the graph and the
corresponding matrix of eigenvectors $X \in \RR^{n \x k}$. The partitioning information is then
computed from  the eigenvectors, using a clustering method for $n$
points in $\RR^k$ (in the literature usually the Kmeans algorithm).

The eigenvalue problem can be formulated,
\begin{equation*}
  \label{eq:minLapl-k-0}
  \min_{Y^T Y = I} \| I - Y^T  A Y \|,
\end{equation*}
($A$ is the adjacency matrix and $\| \cdot \|$ denotes  the Frobenius norm),
which is solved by the eigenvector matrix $X$. We define the minimum in
\eqref{eq:minLapl-k-0}  as a measure of the distance of
the graph to being $k$-partitionable.

The method presented in this paper is  based on  two
observations. Firstly, the information needed for spectral
$k$-partitioning  is contained, not in the eigenvectors as such, but
rather the subspace spanned by the eigenvectors. That subspace can be
represented by $XQ$ for any orthogonal matrix $Q$. Secondly, any
$k$-partitioning can be encoded by a matrix $\Psi \in \RR^{n \x k}$,
where $\psi_{ij} \neq 0$    if node $i$ belongs to partition $j$,
and $\psi_{ij}=0$ otherwise. Thus, given the matrix $X$ of
eigenvectors, we will solve approximately the problem,
\begin{equation}
  \label{eq:minXQPsi-0}
  \min_{Q, \, \Psi} \| XQ - \Psi \|,
\end{equation}
with suitable constraints on $Q$ and $\Psi$. The partitioning is then
obtained directly from $\Psi$. An algorithm for solving
\eqref{eq:minXQPsi-0} was given in \cite{eltr18}. It is the basic
ingredient in the novel fast and reliable multiway partitioning
algorithm. Two cut functions are defined, based on the representation
of the partitioning as the matrix $\Psi$. From these definitions
Cheeger inequalities can be proved, and thereby the relation between the
eigenvalue problem and the partitioning problem is established.

The paper is organized as follows. After a short description of
notation and relevant concepts in Section \ref{sec:notation}, we
briefly review  the classical 2-partitioning problem in Section
\ref{sec:2-part}. The $k$-partitioning problem is  introduced
in Section \ref{sec:multiway}. Our formulation of the partitioning
problem leads naturally to a cut function, for which a Cheeger
inequality follows directly, Section \ref{sec:Psicut}. An alternating
method for solving 
\eqref{eq:minXQPsi-0} is described in Section \ref{sec:SSO}, and the
computation of a starting approximation \cite{dmy19} in Section 
\ref{sec:CPQR}. Another cut function that gives a tighter Cheeger
inequality  is presented in Section
\ref{sec:Phicut}. Numerical examples with data from applications are
given in Section \ref{sec:num-exp}; there variants of the new method
are compared to each other and to Kmeans clustering. In Section
\ref{sec:conclusions} some final  conclusions are  given. 

Surveys of graph partitioning can be found in \cite{lux07,fohr16}.
Early results on multiway spectral partitioning and cut functions are
given in \cite{mexu04}.  The multiway partitioning algorithm described
in this paper uses the one in \cite{dmy19} as a starting
approximation. It is also related to the algorithms in
\cite{zhjo08,pasw22,psz15}, see Section \ref{sec:multiway}.  Cut
functions and Cheeger inequalities for  $k$-partitioning 
are discussed in \cite{lrtv12}. In this paper we assume that the
number of clusters $k$ is given a priori. The problem of choosing $k$
from the data is discussed in \cite{aegl21}.

\section{Notation and Preliminaries}
\label{sec:notation}

The identity matrix is denoted $I$; its dimension will be clear from
the  context. 
The  group of orthogonal matrices  $Q \in \RR^{k \x k}$ is denoted
$\cO(k)$.  A matrix $X \in \RR^{n \x k}$ is said to be in
\emph{indicator form}\footnote{Our definition is less restricted
  than   that in e.g. \cite{psz15,dmy19}, where the nonzero elements
  of the columns   are assumed to be equal.},   if, for some 
permutation matrix $P$,  we can write it as
\[
  X =
  P  \begin{pmatrix}
    x_1 & 0 & \cdots & 0\\
    0 & x_2 & \cdots & 0\\
    \vdots & & \ddots & \vdots\\
    0 & 0 &\cdots & x_k
  \end{pmatrix},
\]
where all components of the vectors $x_i$ are nonzero.
The set of matrices  $X \in \RR^{n \x k}$, which are in indicator form
and satisfy $X^TX = I$, is denoted $\cOI(n,k)$. If 
$X \in \cOI(n,k)$ is nonnegative, then we say $X \in
\cOIplus(n,k)$.

A $k$-partitioning of the integers $1,2,\ldots,n$ is a grouping
into $k$ disjoint and non-empty subsets, such that the union of the
subsets is equal to the whole set. Such a grouping can be encoded as
a matrix in $\cOI(n,k)$.  

We will assume that all eigenvectors of symmetric matrices are
normalized to Euclidean length 1. We  define unnormalized constant vectors, 
\begin{equation}
  \label{eq:gi}
  g_i= 
  \begin{pmatrix}
  1\\
  1 \\
  \vdots \\
  1
\end{pmatrix} \in \RR^{n_i},  \qquad
 g= 
  \begin{pmatrix}
  1\\
  1 \\
  \vdots \\
  1
  \end{pmatrix} \in \RR^{n}.  
\end{equation}
 Throughout  we will use the Frobenius matrix norm, $\| A
 \| = (\sum_{i,j} a_{ij}^2)^{1/2}$, and corresponding inner product
 $\langle A , B \rangle = \tr(A^T B)$.  

Let the symmetric, nonnegative matrix 
$B \in \RR^{n \x n}$ be the \emph{unnormalized adjacency matrix} of an undirected graph
$\cG$ over $n$ nodes (vertices). This means that an element $b_{ij}$ 
 is nonnegative, if there is an edge between nodes $i$
and  $j$, otherwise $b_{ij}=0$.   We let 
 $b_{ij}$ be  the   edge weight. Throughout we assume that the diagonal of
the adjacency matrix is zero, i.e., no self-loops are allowed in the
graph. 

 Let $d=Bg$ be the degree vector, where $d_i$
is the degree, i.e. the number of edges, of node $i$.  The volume of the graph is defined
$\Vol(\cG)=\sum_{i=1}^n d_i$.  Define the diagonal matrix $D =
\diag(d)$; thus  we have $B g = D g$.    The \emph{normalized
  adjacency matrix} is defined
\[
  A = D^{-1/2} B D^{-1/2}.
\]
Denote the eigenvalues of $A$ by $\lambda_i, \, i=1,2,\ldots,n$, and
assume the eigenvalues are ordered,
\[
  \lambda_1 \geq \lambda_2 \geq \cdots \geq \lambda_n.
\]
The normalized Laplacian matrix for $\cG$ is $I - A$. Denote its
eigenvalues by $\mu_i$.   The eigenvalue problems for  the normalized
adjacency and the normalized Laplacian matrices  are
 equivalent. The eigenvectors are the same, and we have 
 $\lambda_i=1-\mu_i, \, i=1,2,\ldots,n$. In this paper we will
use the normalized adjacency matrix\footnote{For large and sparse
  adjacency matrices it usually easier to compute the
  largest eigenvalues  of the adjacency matrix than
  the smallest eigenvalues  of the Laplacian.}; the equivalence
to the Laplacian formulation gives automatically the connection to the
cost (cut)  functions for partitioning the graph.

A subgraph is a subset of the nodes and corresponding edges.  The
graph is called connected if there is no subgraph isolated from the
rest of the graph.  Connectedness is equivalent to irreducibility
of the adjacency matrix. 
Any symmetric matrix $A$ is called
\emph{reducible},
if there exists a permutation matrix $P$ such that $P A P^T$ is
block-diagonal,
\[
  P A P^T =
  \begin{pmatrix}
    A_1 & 0 \\
    0 & A_2
  \end{pmatrix}.
\]
If there is no such permutation, then the matrix is called
\emph{irreducible}. Thus an undirected graph is connected if and only if its
adjacency matrix is irreducible. If the adjacency matrix is reducible,
and can be permuted to block-diagonal form with $k$ square blocks,
 we will call it \emph{$k$-reducible}. We will use the corresponding term
\emph{$k$-partitionable} for  a graph that  has at least $k$
components.
Throughout we will identify the graph with its adjacency matrix, and
choose  the concepts and notation that are
most convenient in the actual situation.

We now recall some basic properties of the eigenvalue problem for
the normalized adjacency matrix $A$ and its relation to connectedness
of the graph (for an elementary introduction,
see e.g. \cite[Chapter 10]{eldenmatrix19}, and for a  comprehensive
treatment, see \cite[Chapter 8]{mey00}, \cite{lux07}).

The largest eigenvalue of $A$ is $\lambda_1=1$, with eigenvector
$x_1=\alpha D^{1/2} g$ ($\alpha$ is a normalization constant). Clearly $x_1$   is positive (it is called the Perron
vector). If $A$ is irreducible 
($\cG$ is connected), then $\lambda_1$ is a simple eigenvalue, and
$x_1$ is the only positive eigenvector. The second eigenvalue is
strictly less than 1, and, due to orthogonality, $x_1^T x_2 = 0$,
and $x_2$ (the Fiedler vector \cite{fied73})  must have positive and negative
components.


\section{Spectral 2-partitioning}
\label{sec:2-part}

The purpose of spectral 2-partitioning is to find a partitioning of
the the nodes in the graph so that the cost for splitting the graph in
two is as small as possible. As  cut  function it is common to use
\emph{conductance}. Let $\cR$ be a set of nodes and $\cR^C$ its
complement.  Then the conductance of the partitioning is defined
\cite{chung97}, 
\[
  h(\cR) = \frac{| E(\cR,\cR^C) | }{\min(\Vol(\cR),\Vol(\cR^C))} \, ,
\]
where $| E(\cR,\cR^C) |$ is the sum of the weights of all edges that
start in $\cR$ and end in $\cR^C$. The \emph{Cheeger constant} of the
graph is
\[
  h_\cG = \min h(\cR),
\]
where the minimization is done over all partitionings. 
Spectral 2-partitioning is based on the fact that the closeness of the
graph to disconnectedness is related to the smallness of the quantity
$1-\lambda_2$. This is expressed in the \emph{Cheeger inequality}
\cite{chung97},
\begin{equation}
  \label{eq:cheeger-2}
  \frac{h_\cG^2}{2}  < 1 - \lambda_2 \leq 2 h_\cG.
\end{equation}
We can interpret this as follows:
\begin{quote}\emph{
  There is a partitioning of the graph with small conductance if and
  only if $1-\lambda_2$ is small.}
\end{quote}

Furthermore, an algorithm for computing a partitioning that is close
to optimal is based on properties of the eigenvectors: $x_1>0$ and
$x_1^T x_2 =0$.    The algorithm is sketched in Algorithm 1.

\begin{algorithm}[htb]\label{alg:2-part}
{\vskip 4pt}
\begin{algorithmic}
{\vskip 4pt}
  \STATE{\textbf{1.} Compute the two largest eigenvalues  and
    corresponding eigenvectors $x_1$ and $x_2$ of $A$.  }
  \STATE{\textbf{2.} Reorder the elements of the Fiedler vector $x_2$
    in ascending 
    order. Apply the reordering to the adjacency matrix. }
  \STATE{\textbf{3.} Compute the conductance for partitioning the graph in
    the neighborhood of the sign change of the reordered vector, and split
    where the conductance is smallest. }
\end{algorithmic}

\caption{\textbf{Spectral 2-partitioning:} Given a normalized
  adjacency matrix $A \in \RR^{n \x k}$, compute a reordering and
  partitioning, such that the cost for splitting the graph is low.} 
\end{algorithm}

Assume for the moment that the graph is disconnected, having two
components, each of which is well connected. The  normalized
adjacency matrix can be written
\[
  A =
  \begin{pmatrix}
    A_1 & 0 \\
    0 & A_2
  \end{pmatrix}.
\]
A has a double eigenvalue 1, and the nonunique eigenvectors can be
written in indicator form
\[
  X =
  \begin{pmatrix}
    x_1 & 0 \\
    0 & x_2
  \end{pmatrix}, \quad x_i >0, \; i=1,2,
\]
where $x_1$ and $x_2$ are the Perron vectors of $A_1$ and $A_2$,
respectively.  
If we add an edge with small weight $\varepsilon$ between the components, 
then, from the perturbation theory for the symmetric eigenvalue
problem \cite[Chapter 8.1]{govl13}, the corresponding adjacency matrix
$A_\varepsilon$  has the 
eigenvalues 1 and $\lambda_2 \approx 1 - \varepsilon <1$, and  the unique
eigenvectors are given by 
\[
  X_\varepsilon    \approx \frac{1}{\sqrt{2}}
  \begin{pmatrix}
   x_1  & -x_1 \\
    x_2 & x_2
  \end{pmatrix}.
\]
The important observation is that while the eigenvectors have been
rotated by approximately $45^{\circ}$ from indicator form, the subspace
spanned by the eigenvectors has changed by $\varepsilon$ \cite[Chapter
8.1]{govl13}.  
So instead
of using the Fiedler vector in Algorithm 1, we could try to determine
a  rotation $Q$ that will make $X_\varepsilon Q$ close to indicator
form, and obtain the partitioning information from that. This idea
can be generalized to the case $k>2$.

\section{ Multiway Spectral Partitioning}
\label{sec:multiway}

To find more than two subgraphs, the 2-partitioning algorithm can be
applied recursively. However, this is reported to give poor solutions
in some common cases \cite{rine14}.  Therefore  it is often preferable to
try to find a 
partitioning into $k>2$ subgraphs at the same time. A natural
generalization of 2-partitioning is to compute the $k$ largest
eigenvalues of $A$,
and the corresponding matrix of eigenvectors
$X = (x_1 \, x_2 \, \ldots \, x_k) \in \RR^{n \x k}$. However, the
useful sign pattern of the first two eigenvectors does does not carry
over  in a
 way that can be easily used to find the $k$-partitioning. A standard
method is to apply the Kmeans algorithm to the rows of $X$ (see
e.g. \cite{dhil01a,lux07}). However, Kmeans has the
disadvantage that  different applications of the algorithm  to the same
data often give different results, cf. Section \ref{sec:num-exp}. 

We will now describe the ideas behind spectral partitioning for $k>2$,
and we will emphasize the well-known invariance properties of the underlying
eigenvalue problem.
 Consider the unnormalized adjacency matrix $ B$  of a connected graph on
 $n$ nodes, and  define the
normalized adjacency matrix
\begin{equation}
  \label{eq:Atk}
   A =  D^{-1/2}  B  D^{-1/2};
\end{equation}
 the normalization matrix $ D$ is the standard degree matrix,
 satisfying $ D g =  B g$. 
With this normalization, the largest eigenvalue of $ A$ is 
equal to 1, and, due to the connectedness of the graph,
\[
  \lambda_\nu < 1, \quad \nu=2,3,\ldots,k.
\]
Throughout we will assume that $\lambda_k > \lambda_{k+1}$.

The basis of spectral partitioning for $k>2$ is to compute
the $k$ largest eigenvalues of $ A$,  and the corresponding eigenvectors. 
The latter can be found by considering matrices $Y \in \RR^{n \x k}$
and the minimization problem 
\begin{equation}
  \label{eq:minLapl-k}
  \min_{Y^T Y = I} \| I - Y^T  A Y \|.
\end{equation}
It is well-known that the eigenvector matrix $X$ solves this problem,
and the minimum is equal to
\begin{equation}
  \label{eq:Lk}
  \left(\sum_{\nu=2}^k (1 - \lambda_\nu)^2 \right)^{1/2} =: L_k,
\end{equation}
(since $\lambda_1=1$ ).  Clearly, the solution is not unique: any
matrix $Y = X Q$, for arbitrary $Q \in \cO(k)$ is a minimizer.  If
instead we define the minimization problem \eqref{eq:minLapl-k} over
subspaces in $\RR^n$ of dimension $k$ (the Grassmann manifold
\cite{eas98}), then, under the assumption 
$\lambda_k > \lambda_{k+1}$, the problem has a unique solution,
which is the subspace spanned by $X$.

The following simple result is a direct consequence of the basic
properties of the eigenvalue and partitioning problems.

\begin{proposition}
  \label{prop:Lk}
  Let $ A$ be the normalized adjacency matrix of an undirected
  graph, and assume that $\lambda_k > \lambda_{k+1}$. Then $L_k=0$  if
  and only if $ A$ is $k$-reducible (the graph is
  $k$-partitionable).   
\end{proposition}


From the perturbation theory of the symmetric eigenvalue problem
\cite[Chapter 8.1]{govl13} the eigenvalues
$\lambda_1,\ldots,\lambda_k$ are continuous functions of perturbations
of the adjacency matrix (as long as $\lambda_k > \lambda_{k+1}$).
Therefore, in analogy to Fiedler's definition of $(1-\lambda_2)$ as
algebraic connectivity \cite{fied73} for $2$-partitioning, we can use $L_k$ as a measure
of the distance of a graph to being $k$-partitionable. See also
Appendix \ref{app:Lk} for an illustration of the properties of $L_k$. 

Clearly, $L_k$  does not depend on the ordering of the nodes;
equivalently, it does not depend on symmetric permutations of $B$. We
will see that the same is valid for the cut functions that we will
define. Therefore, it is no restriction to discuss and illustrate the
partitioning problem in terms of any particular ordering.  We write
 \begin{equation}
   \label{eq:BDBE}
    B = D_B + E,
 \end{equation}
where
\begin{equation}
  \label{eq:B0k}
  D_B =
  \begin{pmatrix}
    B_1 & 0 & \cdots & 0 \\
    0 & B_2 &        & 0 \\
    \vdots & & \ddots & \vdots \\
    0      & 0  & \cdots & B_k
  \end{pmatrix}, \qquad 
  E =     \begin{pmatrix}
    0 & E_{12} & \cdots & E_{1k} \\
    E_{12}^T & 0 &        & E_{2k} \\
    \vdots & & \ddots & \vdots \\
    E_{1k}^T& E_{2k}^T  & \cdots & 0
  \end{pmatrix}.
\end{equation}
We assume that for $i=1,2,\ldots,k,$  $B_i \in \RR^{n_i \x n_i}$,
with $\sum_{i=1}^k n_i =n.$ Note that $B$ is irreducible (and the
corresponding graph is connected) if at least one matrix $E_{ij}$ or
$E_{ij}^T$ in each block row is nonzero.

The  spectral partitioning problem can be loosely formulated:
\begin{quote}
\emph{Find a permutation $P$ for some original adjacency matrix $\tilde
  B$,  giving $B=P \tilde B P^T$, and a 
 $k$-partitioning, such that,  in some sense 
 (which will soon be clear), $D_B$ is considerably larger than $E$.}
\end{quote}


The partitioning of the graph, and the equivalent blocking of the
adjacency matrix, can be represented in terms of a matrix $\Psi \in
\RR^{n \x k}$ in indicator form, where $\psi_{ij} \neq 0$, if node $i$
belongs to partition $j$, and $\psi_{ij}=0$ otherwise.
This indicator matrix  will  be computed from the eigenvector
matrix $X$ that solves \eqref{eq:minLapl-k}. But first we will derive
a cut  (cost) function for the partitioning of the graph. We will come
back to cut functions in Section \ref{sec:Phicut}; here it is enough to
say that the cut function is a measure of how much the partitioning
defined by \eqref{eq:BDBE}-\eqref{eq:B0k} deviates from
block-diagonal. 

\subsection{A Cut Function}
\label{sec:Psicut}

As indicator matrices are crucial in the partitioning problem, we will evaluate
$\| I - Y^T A Y \|$ for one particular indicator matrix $Y$, satisfying
the constraint $Y^T Y = I$. First write the matrix $D$ using the same partitioning
as in \eqref{eq:B0k},
\[
    D =
  \begin{pmatrix}
    D_1 & 0 & \cdots & 0 \\
    0 & D_2 &        & 0 \\
    \vdots & & \ddots & \vdots \\
    0      & 0  & \cdots & D_k
  \end{pmatrix},
\]
put
\[
  \omega_i = g_i^T D_i g_i,  \quad i=1,2,\ldots,k, 
\]
(recall \eqref{eq:gi}),  and define the matrices
\begin{equation}
  \label{eq:OmegaG}
  \Omega=
  \begin{pmatrix}
    \omega_1 I & 0 & \cdots & 0 \\
    0 & \omega_2 I &        & 0 \\
    \vdots & & \ddots & \vdots \\
    0      & 0  & \cdots & \omega_k I
  \end{pmatrix} \in R^{n \x n},
  \quad
  G = 
  \begin{pmatrix}
    g_1 & 0 & \cdots & 0 \\
    0 & g_2 &        & 0 \\
    \vdots & & \ddots & \vdots \\
    0      & 0  & \cdots & g_k
  \end{pmatrix} \in \RR^{n \x k}.    
\end{equation}
Then put
\[
  Y = D^{1/2} \Omega^{-1/2} G.
\]
Clearly $Y$ is positive and  in indicator form.  It is straightforward
to verify that $Y^T Y = I$.  Using \eqref{eq:Atk} and \eqref{eq:BDBE} we
can write
\begin{equation*}
  I - Y^T A Y  = I - G^T \Omega^{-1/2} D_B \Omega^{-1/2} G -
                G^T \Omega^{-1/2} E \Omega^{-1/2} G =: \Gamma,
\end{equation*}
where the symmetric matrix $\Gamma \in \RR^{k \x k}$ is given by
\begin{equation}
  \label{eq:gammaij}
  \gamma_{ij} =
  \begin{cases}
    1 - {g_i^T B_i g_i}/{\omega_i}, & i=j,\\
      - {g_i^T E_{ij} g_j}/{(\omega_i \omega_j)^{1/2}}, & i \neq j.
  \end{cases}
\end{equation}
We  now define a cut function for a given $k$-partitioning
\eqref{eq:B0k}, 
\begin{equation}
   \label{eq:PsiCut}
      \PsiCut = \| \Gamma \|. 
 \end{equation}
Note that $\gamma_{ii}$ is equal to 1 minus the sum of the weights of
edges within partition  $i$, divided by the total weights of edges within
and out from partition $i$. So $\gamma_{ii} \geq 0$. Similarly,
$\gamma_{ij}$ for $i \neq j$  is minus the sum of the weights of the edges between
partitions $i$ and $j$, divided by the square roots of the total
weights of the partitions. This description of $\PsiCut$ shows
that the function is invariant under renumberings of the
nodes. Therefore it is no restriction to define the function in terms
of the blocking of the adjacency matrix for a particular
numbering.

Further define  
\begin{equation}
  \label{eq:PsiCutG}
\PsiCut_\cG =   \min \PsiCut,
\end{equation}
where the minimum is taken over all $k$-partitionings of  the
graph. Obviously it is NP-hard to compute $\PsiCut_\cG$. 
It is easy to see that  $\cG$ is $k$-partitionable if and only if
$\PsiCut_\cG=0.$
Clearly we have the \emph{Cheeger inequality},
\begin{equation}
  \label{eq:cheeger-k}
   L_k = \left(\sum_{\nu=2}^k (1 - \lambda_\nu)^2 \right)^{1/2}
   \leq \PsiCut_\cG.
 \end{equation}
 We can interpret the inequality as follows:
 \begin{quote} 
  \emph{If there exists a partitioning for which $\PsiCut$ is small,
     then the $k$ largest eigenvalues of the normalized adjacency
     matrix  $A$ are close to 1.}
 \end{quote}

 In the  graph partitioning literature a function of the
 type of $\PsiCut$ is  often referred to as a \emph{cut}, see
 e.g. \cite{mexu04,zhjo08,lgt14}, or \emph{expansion constant} \cite{psz15}. 
Our cut function $\PsiCut$ is very natural to the problem, as the Cheeger
inequality \eqref{eq:cheeger-k} follows automatically from  the
definition. Another related cut function will be given in Section
\ref{sec:Phicut}. 

We will now present a partitioning algorithm.
The key  observation is that any spectral partitioning algorithm is
equivalent to computing an indicator matrix from the eigenvector
matrix $X$.    At first sight, the eigenvector
matrix $X$ is far from being in indicator form (e.g. the first eigenvector
is positive).  But the Cheeger inequality \eqref{eq:cheeger-k}
does not depend on which matrix is used to represent the subspace that is
the solution of \eqref{eq:minLapl-k}. Instead of $X$ we may
consider $XQ$ for some orthogonal matrix $Q$. Thus we may ask: \emph{Is there a
  $Q$ such   that $XQ$ is   close to indicator form?}

We know that if $A$ is $k$-reducible, then the eigenvector matrix can be
written in indicator form. If $\lambda_k-\lambda_{k+1}>0$, and this difference is
not small, then a small
perturbation of $A$ gives a small perturbation of the subspace spanned
by the eigenvectors \cite[Chapter 8.1.3]{govl13}. 
Therefore,   if
$A$ is close to $k$-reducible, then there will exist a $Q$  such
that $XQ$ is close to indicator form. 
So we are lead to the minimization problem,
\begin{equation}
  \label{eq:minXPsiQ}
  \min_{Q,\Psi} \| X Q - \Psi \| = \min_{Q,\Psi} \| X - \Psi Q^T \|,
  \qquad Q \in \cO(k),
\end{equation}
where $\Psi$ is required to be in indicator form.  If we can solve
\eqref{eq:minXPsiQ}, then the 
partitioning can be read directly from the non-zeros of $\Psi$.
In the following
section we describe  Algorithm 3, for solving \eqref{eq:minXPsiQ}.
Using that, we state a simple $k$-partitioning algorithm in Algorithm
2.

\begin{algorithm}[htb]\label{alg:k-part}
{\vskip 4pt}
\begin{algorithmic}
{\vskip 4pt}
  \STATE{\textbf{1.} Compute  the eigenvector matrix $X \in \RR^{n \x k}$
       corresponding to the $k$ largest eigenvalues of $A$. }
  \STATE{\textbf{2.} Solve \eqref{eq:minXPsiQ}  using Algorithm 3.} 
  \STATE{\textbf{3.} Find the partitioning from $\Psi$, and  compute
    the cut function $\PsiCut$.} 
\end{algorithmic}
\caption{\textbf{ $k$-partitioning algorithm:} Given the normalized
  adjacency matrix $A$ of a graph, compute a $k$-partitioning. }
\end{algorithm}

\subsection{A Semi-Sparse
  Orthogonal    Approximation of $X$. }
\label{sec:SSO}

Consider the optimization problem  \eqref{eq:minXPsiQ}.  
 We want to approximate the  matrix of eigenvectors by a product of
 matrices, where  the first factor is in indicator form, and the
 second is an orthogonal 
matrix. This can be called a \emph{semi-sparse 
  orthogonal approximation}.
To our knowledge there is no explicit solution of this problem.
In \cite{eltr18} we give an alternating algorithm.  Assume that $Q$ is
given,   put $Y=X Q$ and consider  the problem 
\begin{equation}
  \label{eq:Psicomp-2}
  \min_\Psi \| Y - \Psi \|,
\end{equation}
where $\Psi$ is to be  in indicator form. This  constraint is
equivalent to the 
requirement that exactly one element of each row of $\Psi$ is
nonzero. Therefore \eqref{eq:Psicomp-2}
consists of $n$ independent minimization problems, one for each row,
and the solution is given by 
\begin{equation}
  \label{eq:Psicomp-3}
  \Psi_{ij}=
  \begin{cases}
    y_{i \nu_i}, & j=\nu_i,\\
    0, & \text{otherwise.}
  \end{cases}, \qquad i=1,2,\ldots,n,
\end{equation}
where
\begin{equation}
  \label{eq:Psicomp-1}
  \nu_i = \arg\max_\nu | y_{i\nu} |, \qquad i=1,2,\ldots,n.  
\end{equation}
Then assume that $\Psi$ is given. The minimization  problem 
\begin{equation}
  \label{eq:procrustes}
  \min_{Q}\| X  - \Psi Q^T \|, \qquad Q \in \cO(k),
\end{equation}
is an orthogonal Procrustes problem \cite[Section
6.4.1]{govl13}. Let $\Psi^T X = U \Sigma V^T$ be the Singular Value
Decomposition  (SVD). The solution of \eqref{eq:procrustes}
is $Q=V U^T$. 

Starting with an initial approximation of $\Psi$ the algorithm
alternates between solving the problems for $Q$ and $\Psi$, and stops
when the gradient of the objective function 
is small enough in norm. As is common with alternating algorithms, the
convergence  rate is linear. It is of course crucial to have a good
initial approximation; we will come back to this in the following
subsection.

\begin{algorithm}[htb]\label{alg:SSO}
{\vskip 4pt}
\begin{algorithmic}
{\vskip 4pt}
  \REPEAT
  \STATE{\textbf{1.} Compute $Q$ as the solution of  \eqref{eq:procrustes}.}
      \STATE{\textbf{2.} Put $Y=XQ$ and compute $\Psi$ from
        \eqref{eq:Psicomp-2}-\eqref{eq:Psicomp-1}. }
  \UNTIL{convergence}      
\end{algorithmic}

\caption{\textbf{SSO algorithm:} Given   $X \in
  \RR^{n \x k}$ and  an initial approximation $\Psi$ in indicator
  form, compute an SSO 
  approximation $X \approx \Psi Q^T$.} 
\end{algorithm}

As the optimization problem is non-convex, we cannot guarantee that
the computed solution is globally optimal.  To check if we are close
to a local minimum, we  compute gradients. Let
\[
  r(\Psi,Q) = \frac{1}{2} \| X- \Psi Q^T \|^2,
\]
and denote the gradients $\nabla_\Psi$ and $\nabla_Q$ (the
$Q$-gradient must take into account that $Q \in \cO(k)$).   

\begin{proposition}
  \label{prop:gradients}
  Assume that the indicator structure of $\Psi$ is fixed, and let
  $\Psi$ and $Q$ be the output of Algorithm 3.
  Then
  \begin{equation*}\label{eq:gradients}
    \nabla_\Psi = 0, \qquad \| \nabla_Q \| = \frac{1}{2}\| \Psi^T XQ -
    (XQ)^T \Psi \|. 
  \end{equation*}
\end{proposition}

The proof is given in Appendix \ref{app:gradients}.

To our knowledge the idea to use an alternating algorithm for computing the
partitioning was first used in \cite[Algorithm 1]{zhjo08} (see also
\cite[Section 3.4]{pasw22}). However, there it was assumed that the indicator
matrix had the form $\Psi = P G$, with a permutation matrix $P$,
and $G$ as in \eqref{eq:OmegaG}. 

The approximation can also be written
\begin{equation}
  \label{eq:XapprSQ}
  X \approx \Phi S Q^T,
\end{equation}
where $\Phi \in \cOI(n,k)$,  $S$  is diagonal,
$S=\diag(s_1, s_2, \ldots, s_k)$. We use the convention that the $s_i$
are positive and ordered by magnitude. 
The following lemma will later be needed. 

\begin{lemma}
  \label{lem:S<=1}
  The quantities $s_i$ satisfy
  \[
    1 \geq s_1 \geq s_2 \geq \cdots \geq s_k.
  \]
\end{lemma}

\begin{proof}
  The matrix $Y$ in \eqref{eq:Psicomp-2} has columns of length 1. The
  columns of $\Psi$ are obtained by putting some of the elements of
  $Y$ equal to zero. Therefore the length of column $i$ of $\Psi$,
  which is the the quantity $s_i$, is less than or equal to 1.
\end{proof}

\subsection{A   Starting Approximation }
\label{sec:CPQR}

In \cite{dmy19} a direct multiway spectral partitioning method, called
CPQR, is described.  It computes a  QR decomposition with column pivoting
\cite[Section 5.4.2]{govl13} and the polar decomposition \cite[Section
9.4.3]{govl13}. The method is summarized in Algorithm 4. 

\begin{algorithm}[htb]\label{alg:CPQR}
{\vskip 4pt}
\begin{algorithmic}
{\vskip 4pt}
  \STATE{\textbf{1.} Compute the QR decomposition with column pivoting
  \[
    X^T \Pi_0 = U R,
  \]
  where $\Pi$ is a permutation matrix.}
      \STATE{\textbf{2.} Partition
  \[
    X^T \Pi_0 = (Z_1 \; Z_2), \quad Z_1 \in \RR^{k \x k},
  \]
  and compute the SVD,  $Z_1 = U_Z \Sigma_Z V_Z^T$, giving  the polar
  decomposition 
  \[
    Z_1 = Q_0 H, \qquad Q_0 = U_Z V_Z^T, \quad H = V_Z \Sigma_Z V_Z^T,
  \]
  where   $Q_0$ is orthogonal and $H$ is symmetric.} 
  \STATE{\textbf{3.} Compute $X Q_0$ and determine a matrix $\Psi_0$
    in indicator 
    form as in     \eqref{eq:Psicomp-1}-\eqref{eq:Psicomp-3}. }
  \end{algorithmic}
\caption{\textbf{CPQR algorithm} \cite{dmy19}: Given the   matrix of
  eigenvectors $X \in   \RR^{n \x k}$,  compute a semi-sparse
  approximation  $X \approx \Psi_0 Q_0^T$.   } 
\end{algorithm}

 The QR decomposition is only used to find the $k$ rows of $X$
that ``dominate'', in the sense that they are close to being the best
basis vectors in  $\RR^k$ among all row vectors (cf. the discussion in
\cite[Section 5.4]{govl13}). Therefore, these vectors should be good
for determining an  orthogonal matrix $Q_0$ such that $X Q_0$ is close
to indicator form, and the corresponding $\Psi_0$ should be a good 
starting point for a partitioning algorithm. This is confirmed in our
numerical experiments. It is suggested in
\cite{dmy19} to use CPQR as a  starting point for spectral
partitioning by Kmeans
clustering.  We use  $Q_0$  as 
starting approximation for SSO.

In \cite{dmy19} it is proved that if $X$ is close to an indicator
matrix $W$ (in the restricted sense of this term, i.e., the nonzero
elements of each column are equal), then there is a
permutation matrix $\Pi$ such that $X Q_0$  is  close to $W \Pi$.   
Therefore it is not surprising that it also gives a good
approximation, 
\begin{equation}
  \label{eq:XapprPsiQt}
  X \approx \Psi_0 Q_0^T = \Phi_0 S_0 Q_0^T,  
\end{equation}
where  $\Psi_0$ and $\Phi_0$  are   indicator matrices in the sense of
this paper,  $\Phi_0 \in \cOI(n,k)$, and $S_0 \in \RR^{k \x k}$ is
diagonal (as in \eqref{eq:XapprSQ}). 

If we apply the CPQR algorithm to a matrix of eigenvectors, which
satisfies \eqref{eq:XapprPsiQt} with equality, then, naturally, it recovers
the exact solution. This follows from the results in
\cite{dmy19}; it is also straightforward  to give a a constructive
proof. 
Therefore, for a general disconnected graph, CPQR can be used to
compute the indicator matrix for partitioning the graph into its
components\footnote{The problem can alternatively be solved by an
  algorithm of Tarjan \cite{tarj72}. A Matlab implementation by
  D. Gleich is available at
  \url{https://se.mathworks.com/matlabcentral/fileexchange/
    24134-gaimc-graph-algorithms-in-matlab-code}.}.

\subsection{Another Cut Function}
\label{sec:Phicut}
As we noted in the introduction, cut functions are often used in the
literature, first as a motivation for the spectral partitioning
approach, and then as a measure to evaluate the quality of a computed
$k$-partitioning. In this paper we emphasize the second use of the
function.

Let $\Theta$ denote a function, defined  for any 
given $k$-partitioning of a graph $\cG$. Let 
\[
  \Theta_\cG  = \min \Theta,
\]
where the minimum is taken over all $k$-partitionings.

We will call $\Theta$  a
\emph{cut function} if
\begin{enumerate}
\item $\Theta \geq 0$, and it is  invariant under  renumberings of the nodes of the graph,
\item%
  $\Theta_\cG = 0$,  if and only if $\cG$
    is $k$-partitionable.
\end{enumerate}
 It is easy to see that $\PsiCut$ defined in Section
 \ref{sec:Psicut}  is a  cut function.
  
Due to the permutation invariance we can write the cut function in
terms of a  normalized adjacency matrix, in blocked form, 
\begin{equation}\label{eq:Ablocked}
  A = D^{-1/2} B D^{-1/2} =
     \begin{pmatrix}
    A_1 & F_{12} & \cdots & F_{1k} \\
    F_{12}^T & A_2 &  \cdots      & F_{2k} \\
    \vdots & \vdots & \ddots & \vdots \\
    F_{1k}^T& F_{2k}^T  & \cdots & A_k
  \end{pmatrix},
\end{equation}
where the normalization matrix is
\[
  D=
   \begin{pmatrix}
    D_1 & 0  & \cdots & 0  \\
    0  & D_2 &  \cdots      & 0 \\
    \vdots &\vdots & \ddots & \vdots \\
    0 & 0   & \cdots & D_k
  \end{pmatrix}.
\]

For the CPQR and SSO methods we can define a cut function based on the
approximation
\[
  X \approx \Phi S Q^T,
\]
where $\Phi \in \cOI(n,k)$. If $A$
is not far from being $k$-reducible we can expect that only a small
fraction of   the nonzero elements of $\Phi$ are negative (cf. our
numerical experiments in Section \ref{sec:num-exp}). In any case,
define $\Phi_+ = | \Phi |$, where $| \Phi |$  denotes element-wise
absolute value.   Clearly, $\Phi^+ \in \cOIplus(n,k)$.  Define
\begin{equation}
  \label{eq:PhiCut}
  \PhiCut = \| I - \Phi_+^T A \Phi_+ \|. 
\end{equation}
Clearly,   $\PhiCut$ is invariant under symmetric permutations.
The minimal $\PhiCut$ is defined
\begin{equation}
  \label{eq:PhiCut-min}
  \PhiCut_\cG = \min_{\Phi \in \cOIplus(n,k)} \| I - \Phi^T A \Phi
  \|. 
\end{equation}

\begin{proposition}\label{prop:PhiCut}
  $\PhiCut$ is a cut function, and the Cheeger inequality,
  \begin{equation}
    \label{eq:Cheeger-1}
  L_k = \min_ {Y^T Y = I} \| I - Y^T A Y \| \leq \PhiCut_\cG,    
  \end{equation}
is satisfied.
\end{proposition}

\begin{proof}
  We will show that $\PhiCut_\cG=0$ if and only if $\cG$ is
  $k$-partitionable.   Let
  \[
    \Phi_+ =
    \begin{pmatrix}
      \varphi_1 & 0 & \cdots & 0\\
      0    &\varphi_2 & \cdots & 0\\
      \vdots &     & \ddots & \vdots\\
      0 & 0 & \cdots & \varphi_k
    \end{pmatrix} \in \cOIplus(n,k),
  \]
  where the partitioning corresponds to that  in $A$
  \eqref{eq:Ablocked}. Put
  \[
    \Gamma_\Phi = I - \Phi_+^T A \Phi_+.
  \]
  The off-diagonal elements of $\Gamma_\Phi$ are
  \[
    - \varphi_i^T F_{ij} \varphi_j,
  \]
  which is a weighted sum of the elements in $F_{ij}$, so it is zero
  if and only if $F_{ij}=0$.  The $i$'th
  diagonal element of $\Gamma_\Phi$ is
  \[
    1- (0 \cdots \varphi_i^T \cdots 0) A
    \begin{pmatrix}
      0 \\
      \vdots\\
      \varphi_i \\
      \vdots\\
      0
    \end{pmatrix}
    = 1 - \varphi_i^T A_i \varphi_i \geq 0,
  \]
  where the inequality follows from the fact that the eigenvalues
  of $A$ are smaller than or equal to 1.

  If $A$ is $k$-reducible then all $F_{ij}$  are zero, and the minimum
  in \eqref{eq:PhiCut-min} is equal to zero and is attained when the
  $\varphi_i$ are the eigenvectors of the $A_i$. If $A$   is
  $k$-irreducible, then at least one $F_{ij}$  is nonzero, and
  $\PhiCut_\cG>0$.  

 The Cheeger inequality follows immediately from the definition of
 $\PhiCut_\cG$. 
\end{proof}
 
An approximate inequality from the left can also be obtained. It is
based on the approximation $X \approx \Phi S Q^T$ (if $XQ$  cannot be
approximated reasonably well by a matrix in indicator form, then
spectral partitioning is not well motivated).

\begin{theorem}
  \label{thm:cheegerPhi}
  Assume that $\lambda_k > \lambda_{k+1}$, and   let $X$ be a solution
  of \eqref{eq:minLapl-k}, and let  $Z =  \Phi S 
  Q^T$ be an SSO approximation with 
  $\Phi \geq 0$. Further,  let $s_{\min}=\min_i(s_i)$ be the smallest diagonal
  element of $S$. Assume that $L_k<1$. Then 
  \begin{equation}
    \label{eq:cheegerPhi}
    \frac{s_{\min}^4}{1 + 4 \| \Delta^T A X\|} \left(\PhiCut\right)^2
    \leq L_k + \cO(\| \Delta\|^2),
  \end{equation}
  where $\Delta= X - Z$.
\end{theorem}

\begin{proof}
  Putting $L(X)=I - X^T A X$, we have
  \begin{align*}
    L(X)&=I - (Z+\Delta)^TA (Z+ \Delta)\\
    &=  I - Z^T A Z - \Delta^T AX - X^T A \Delta + \cO(\| \Delta\|^2),
  \end{align*}
  and
  \begin{align*}
       \langle L(X), L(X) \rangle  &= \langle L(Z), L(Z) \rangle
         -\langle L(Z), \Delta^T A X \rangle
         -\langle L(Z), X^T A \Delta \rangle \\                             
                                   &-  \langle \Delta^TAX,L(X) \rangle
                                     -\langle X^T A \Delta , L(X)
                                     \rangle + \cO(\| \Delta\|^2)\\
    &= \langle L(Z), L(Z) \rangle
         -\langle L(X), \Delta^T A X \rangle
         -\langle L(X), X^T A \Delta \rangle \\                             
                                   &-  \langle \Delta^TAX,L(X) \rangle
                                     -\langle X^T A \Delta , L(X)
                                     \rangle + \cO(\| \Delta\|^2). 
   \end{align*}
  The Cauchy-Schwarz inequality gives
  \[
     | \langle \Delta^TAX,L(X) \rangle | \leq
     \| \Delta^T A X \| \| L(X) \|,
   \]
   so
   \[
     \langle L(X), L(X) \rangle  \geq \langle L(Z), L(Z) \rangle -
     4 \| \Delta^T A X \| \| L(X) \| + \cO(\| \Delta\|^2).
   \]
   Rearranging this we get
   \begin{align}
     \| L(Z) \|^2 &\leq \| L(X) \|^2 + 4 \| \Delta^T A X \| \| L(X) \|
     + \cO(\| \Delta\|^2) \label{eq:L(X)-1}\\
     &\leq \| L(X) \| (1 + 4 \| \Delta^T A X \| ) + \cO(\|
       \Delta\|^2), \label{eq:L(X)-2}
   \end{align}
   where we have used $L_k = \| L(X) \| < 1$.  
  We now bound $\| L(Z) \|$ from below:
  \[
    \| L(Z) \| = \| I - Q S \Phi^T A \Phi S Q^T \|
    = \| I - S \Phi^T A \Phi S  \| =: \| \Gamma_S \|,
  \]
  where the elements of $\Gamma_S$ are given by
  \[
    \gamma_{ij} =
    \begin{cases}
      1 - s_i^2 \varphi_i^T A_i \varphi_i, & i=j,\\
      - s_i s_j \varphi_i^T F_{ij} \varphi_j, & i \neq j.
    \end{cases}
  \]
  Since from Lemma \ref{lem:S<=1} $s_i \leq 1$, we have
  \[
      1 - s_i^2 \varphi_i^T A_i \varphi_i \geq   1 -  \varphi_i^T A_i \varphi_i \geq 
      s_{\min}^2 (1 -  \varphi_i^T A_i \varphi_i).
    \]
    and, 
    \[
      s_i s_j \varphi_i^T F_{ij} \varphi_j \geq
      s_{\min}^2 \varphi_i^T F_{ij} \varphi_j,
    \]
which gives $\| L(Z) \| \geq s^2_{\min} \PhiCut$.  Combining this and
\eqref{eq:L(X)-1}-\eqref{eq:L(X)-2} we have the inequality \eqref{eq:cheegerPhi}. 
  \end{proof}

Note that since $I = (XQ)^T XQ \approx S \Phi^T \Phi S = S^2$,
$s_{\min}$ is not much smaller than 1. Furthermore, if $\| \Delta \|$
is small, then $\| \Delta A X \|$ is small, since $X$ has orthonormal
columns, and $AX=X \Lambda$, where $\Lambda$  is the matrix of
eigenvalues $\lambda_i$, satisfying $\lambda_i \leq 1$.  
  Thus, from the two Cheeger inequalities \eqref{eq:Cheeger-1} and
  \eqref{eq:cheegerPhi}, we have the following statement:
  \begin{quote}\emph{
    There is a partitioning of the graph with small $\PhiCut$ if
    and only if $L_k$ and $\| \Delta \|$ are small. }
\end{quote}

To compute $\PhiCut$ for a given partitioning we require an
approximation $X \approx \Phi S Q^T$. Thus we can use $\PhiCut$  as a
posteriori measure of the quality of the computed partitioning from
CPQR and the SSO method. If
$\PhiCut$ is close to $L_k$ then it is unlikely that there exists
another partitioning substantially closer to optimal.

\section{Numerical Experiments}
\label{sec:num-exp}

The numerical experiments were performed on a standard desktop
computer using Matlab R2019b.
The three algorithms described in this paper were tested on a few
examples. We compared  1) \textrm{CPQR}, 2) \textrm{SSO}, where $\Psi$
is   initialized as the solution  of $\min_\Psi\|X - \Psi \|$, 3)
\textrm{QR+SSO}:  SSO 
initialized with the CPQR solution, and 4) \textrm{Kmeans} with the
standard Matlab initialization (initial cluster centers are chosen
randomly).

As a measure of the quality of a partitioning, we will use $\PsiCut$
defined in  \eqref{eq:gammaij}-\eqref{eq:PsiCut}, and $\PhiCut$.  
For comparison we also give the common measure \textrm{NCut}
(normalized cut)  \cite{shma00,mexu04,pasw22} for a partitioning of the
unnormalized adjacency 
matrix $B$.   Let $\cR_i$ denote the indices belonging
to partition $i$, and $\cR^C_i$ the indices not in partition
$i$. NCut is given
by  
\[
  \mathrm{NCut}=\sum_{i=1}^k
  \frac{| E(\cR_i,\cR_i^C) |}{\omega_i},
\]
where $| E(\cR_i,\cR_i^C) |$ is defined in Section \ref{sec:2-part}.

In the cases when there was an a priori known correct partitioning, we used
the Rand Index \cite{rand71} to evaluate the quality of the computed
partitioning. The index is described in \cite{zhjo08}: Given
  two partitionings $U$ and $V$, ``let $a$ be the number of pairs of
 objects that are in the same set in $U$  and in the same set
 in $V$, and $b$  the number of pairs of objects that are in
 different sets in $U$  and in different sets in $V$. The Rand
 index is given by RI = $(a + b)/ {n \choose 2}$. If RI = 1, the two
 partitions are identical''. We used a Matlab implementation from
 GitHub\footnote{Chris McComb (2022). Adjusted Rand Index
   (\url{https://github.com/cmccomb/rand_index}), GitHub.  Retrieved October
   26, 2022. }. 

 For \textrm{CPQR}, SSO and \textrm{QR+SSO} we give the norm of the residual
 $\| X Q - \Psi\|/\| X \|$.
 In the SSO iterations we  stopped the iterations when no nonzero
 positions in $\Psi$ changed from one iteration to the next, and $\|
 \nabla_Q \| \leq 10^{-5}$.  The execution times for all methods were
 less than $1$  second for 
all examples. 
 

\paragraph{Synthetic Data}
We generated random sparse adjacency matrices of dimension 11,000, with 6
clusters. A parameter determined the closeness to reducibility.
 In Table \ref{tab:synthetic} we give the results
for an  example, where the matrix was relatively close to
reducibility, and the gap between $\lambda_6$ and $\lambda_7$
was quite large.

\begin{table}[htb]
  \centering
  \begin{tabular}{lcccc}
    &    CPQR  &  SSO &   QR+SSO &  Kmeans \\
    \hline
    Rand Index &   1 &  1&  1 &  0.9997\\
    $\| \nabla_Q\|$ &      0.05&  $0.02 \cdot 10^{-5}$ &$0.02 \cdot 10^{-5}$ \\
    Residual &    0.0928 & 0.0894 &  0.0894\\
Iterations  && 3/1 & 1/1 \\
NCut      &   0.657&  0.657 &  0.657 &  0.657\\
    $\PsiCut$       & 0.305 & 0.305 & 0.305  & 0.304\\
    $\PhiCut$ & 0.275 & 0.275 & 0.275 \\
    \hline
  \end{tabular}
  \caption{Synthetic data I. First 7 eigenvalues: 1, 0.923, 0.904,
    0.886,  0.875,  0.837,  0.538.
    $L_k = 0.265$.  ``3/1'' means 3 iterations before the
    nonzero--zero positions in $\Psi$ stabilized, and 1 iteration after.   }
  \label{tab:synthetic}
\end{table}

It is seen that for this well-conditioned example, all methods
performed very well. 
In the second case, shown in Table \ref{tab:synthetic-2}, we let the
matrix deviate further from 
reducibility, which is seen from the eigenvalues and cut values. Still
the partitionings were very close to  correct. 

\begin{table}[htb]
  \centering
  \begin{tabular}{lcccc}
    &      CPQR &  SSO &  QR+SSO &  Kmeans \\
    \hline
    Rand Index&   0.997&  0.998 & 0.998 & 0.997 \\
    $\|\nabla_Q\|$&   0.172&  $0.3\cdot 10^{-5}$ &  $0.02\cdot 10^{-5}$ \\
Residual  &   0.302 & 0.283 & 0.283\\
Iterations &&  4/1 & 2/1\\
NCut       &  2.95 & 2.95  & 2.95 & 2.94 \\
    $\PsiCut$     &  1.33  & 1.33 &  1.33 &  1.33 \\
    $\PhiCut$  & 1.28 & 1.28 & 1.28 \\
    \hline
  \end{tabular}
  \caption{Synthetic data II. First 7 eigenvalues: 1, 0.549,
    0.494, 0.459, 0.437, 0.382, 0.313.  $L_k = 1.20$.}
  \label{tab:synthetic-2}
\end{table}

\paragraph{Astrophysics Collaboration Network}

This data set\footnote{Downloaded from
  \url{http://snap.stanford.edu/data/ca-AstroPh.html} in October 2022.} is a
collaboration network for the arXiv
Astrophysics category. It has several components;  the largest
components  consists of  17903
nodes. The density of the adjacency matrix is $0.012\%$.
The results are given in Table \ref{tab:Astro}. 
\begin{table}[htb]
  \centering
  \begin{tabular}{lcccc}
    &     CPQR  & SSO &  QR+SSO &  Kmeans \\
    \hline
    Residual  &   0.044 &  0.044 & 0.044\\
    $\|\nabla_Q\|$ & 0.0086&  $0.3\cdot 10^{-10}$& $0.3\cdot 10^{-9}$ \\
    Iterations && 2/1 & 1/1\\
    NCut      &   0.522 & 0.522 &  0.522 & 1.39\\
    $\PsiCut$      & 0.268 & 0.268  & 0.268 & 1.22 \\
    $\PhiCut$ & 0.032 & 0.032 & 0.032 \\
    \hline
  \end{tabular}
  \caption{Astrophysics graph. First 7 eigenvalues:
    $1,  0.994, 0.990, 0.984, 0.983, 0.983,     0.980$.  $L_k=0.031$.  
\label{tab:Astro}  }
\end{table}
The problem to choose the number of partitions $k$ is difficult
because the eigenvalues decay very slowly. We chose $k=6$ as in
\cite{dmy19}. Sometimes Kmeans gave results significantly better (in
terms of NCut and $\PsiCut$)  than the other three methods, sometimes
significantly worse. The behavior of CPQR, SSO, and SSO+CPQR was very
consistent.
\begin{table}[htb]
  \centering
  \begin{tabular}{lccc}
       &SSO&QR$+$SSO&Kmeans\\
    \hline
    CPQR      & 1 & 1      &0.980\\
    SSO    &   & 1      & 0.980\\
    QR+SSO  &   &        & 0.980\\
    \hline
  \end{tabular}
  \caption{Astrophysics graph. Pairwise Rand indices. }
  \label{tab:AstroRand}
\end{table}

In order to compare the partitioning performance of the four methods,
we computed  pairwise Rand indices, presented in Table
\ref{tab:AstroRand}. We see that, in spite of the fact that $\PsiCut$
for Kmeans was much higher than for the other methods, the
partitioning in this particular run   
was almost the same.

\paragraph{Power Grid Graph}
The power grid data\footnote{Downloaded
from \url{http://konect.cc/networks/opsahl-powergrid/} in October 2022.}  are described in \cite{wast98,rine14}. The graph contains
information about the power grid of the Western States in the USA. The
edges are power lines between different types of stations. The
adjacency matrix has dimension 4941. In \cite{rine14} it is
partitioned with $k=4$, but it can be partitioned much further at low
cost, which is seen in the very slow decay of the eigenvalues.  We
used $k=14$. The results are given in Table \ref{tab:power}. 
\begin{table}[htb]
  \centering
  \begin{tabular}{lcccc}
    &     CPQR  & SSO &  QR+SSO &  Kmeans \\
    \hline
    Residual &    0.301& 0.290& 0.290\\
    $\| \nabla_Q \|$&  0.253 &  $0.2\cdot 10^{-5}$ & $0.5\cdot 10^{-5}$ \\
    Iterations && 11/2 & 5/1\\
NCut     &    0.254 & 0.243 & 0.246 &  2.71\\
    $\PsiCut$     &  0.085 & 0.082 & 0.083 &  1.47\\
    $\PhiCut$  & 0.0296 & 0.0284 & 0.0287\\
    \hline
  \end{tabular}
  \caption{Power grid graph. First 15 eigenvalues decay very slowly:
    $ 1, 0.9997,\ldots,
    0.9959, 0.9954$. $ L_k = 0.0092$.  
\label{tab:power}  }
\end{table}

With such a small gap between $\lambda_{14}$ and $\lambda_{15}$ we
cannot expect the subspace to be very well determined. Still CPQR, SSO, and 
QR+SSO consistently  gave practically the same results for different
runs with the same parameters. The results for Kmeans, on the other
hand, always differed  between different runs. 
The partitioning given by Kmeans deviated significantly 
 from that given by the other three methods, see Table
\ref{tab:powerrand}, indicating that it is of lower quality. 
\begin{table}[htb]
  \centering
  \begin{tabular}{lccc}
       &SSO&QR$+$SSO&Kmeans\\
    \hline
        CPQR    & 0.992 & 0.992      &0.903\\
    SSO     &   & 0.999      & 0.902\\
    QR+SSO  &   &        & 0.901\\
    \hline
  \end{tabular}
  \caption{Power grid  graph: Pairwise Rand indices. }
  \label{tab:powerrand}
\end{table}

\paragraph{Yeast  Protein Network} 
We tested a  data set describing  protein interactions contained
in yeast\footnote{Downloaded from
\url{http://konect.cc/networks/moreno_propro/} in October 2022.}.  The largest connected
component had 1458 nodes. Again the eigenvalue decay was very slow. We
chose $k=3$ because there was a slightly larger gap between $\lambda_3$
and $\lambda_4$. The results are given in Table \ref{tab:Yeast}.
\begin{table}[htb]
  \centering
  \begin{tabular}{lcccc}
    &     CPQR  & SSO &  QR+SSO &  Kmeans \\
    \hline
    Residual &    0.202 & 0.188&  0.188 \\
    $\| \nabla_Q \|$ &  0.117 &  $0.9\cdot 10^{-6}$& $0.07\cdot 10^{-5}$\\
    Iterations && 3/1 & 2/2\\
NCut    &  0.145 & 0.143 & 0.143 &  1.291\\
    $\PsiCut$    &  0.106 & 0.105 & 0.105  &  1.178\\
    $\PhiCut$    & 0.0201 & 0.0197 & 0.0197\\
    \hline
  \end{tabular}
  \caption{Yeast data. First 4 eigenvalues:
    $1, 0.9917, 0.9848,  0.9831$.
    $L_k=0.0174$.  
\label{tab:Yeast}  }
\end{table}

Again the results of CPQR, SSO, and QR+SSO were stable for different runs,
while those for Kmeans differed. We also ran the data with
$k=12$. There Kmeans failed.  The Rand indices are given in Table
\ref{tab:randyeast}.
\begin{table}[htb]
  \centering
  \begin{tabular}{lcccc}
    &SSO&QR$+$SSO&Kmeans\\
    \hline
    CPQR      & 0.999 & 0.999      &0.631\\
    SSO    &   & 1      & 0.632\\
    QR+SSO  &   &        & 0.632\\
    \hline
  \end{tabular}
  \caption{Yeast data: Pairwise rand indices. }
  \label{tab:randyeast}
\end{table}
Again Kmeans did not perform as consistently as the other methods and
Tables \ref{tab:Yeast} and \ref{tab:randyeast}  indicate that the
partitioning has significantly lower quality than that given by the other methods.

\paragraph{Twitch Social Network}

We downloaded a Twitch Social network from the SNAP
collection\footnote{\url{http://snap.stanford.edu/data/twitch-social-networks.html},
  downloaded in  December 2022.}. The data have been collected and 
described  in 
\cite{ras21}. The Portuguese language network (PTBR) is a
connected undirected graph on 1912 nodes. We attempted to partition it
for different values of $k$, but the results for all methods were
inconsistent, probably due to the fact that there are many very small
subgraphs that are loosely connected to the rest of the
graph. Therefore we replaced the unnormalized adjacency method $B$ by
$0.999 B + 0.001 g g^T$, i.e. we made the graph complete with edges of
small weight (we also removed self-loops). This technique is called
regularization and is described  e.g. in \cite{zhro18}. Now the
results with CPQR, SSO, and QR+SSO stabilized and were the same in
different runs, as expected. The behavior of Kmeans, on the other
hand, remained very inconsistent. The results for $k=4$  are given
Tables \ref{tab:twitch}-\ref{tab:rand-twitch}. 
\begin{table}[htb]
  \centering
  \begin{tabular}{lcccc}
    &     CPQR  & SSO &  QR+SSO &  Kmeans \\
    \hline
    Residual &    0.268 & 0.254&  0.254 \\
    $\| \nabla_Q\|$& 0.148&  $0.1\cdot 10^{-5}$&  $0.4 \cdot 10^{-5}$\\
    Iterations && 6/2 & 3/2\\
NCut    &  1.58 & 1.58 & 1.58 &  1.55\\
    $\PsiCut$    &  0.936 & 0.935 & 0.935  &  0.915\\
    $\PhiCut$  & 0.719  & 0.718 & 0.718\\
    \hline
  \end{tabular}
  \caption{Twitch data. First 5 eigenvalues:
    $1, 0.750, 0.593, 0.539, 0.472$.
    $L_k=0.664$.  
\label{tab:twitch}  }
\end{table}
\begin{table}[htb]
  \centering
  \begin{tabular}{lcccc}
    &SSO&QR$+$SSO&Kmeans\\
    \hline
    CPQR      & 0.931 & 0.931      &0.757\\
    SSO    &   & 1      & 0.707\\
    QR+SSO  &   &        & 0.707\\
    \hline
  \end{tabular}
  \caption{Twitch data: Pairwise rand indices. }
  \label{tab:rand-twitch}
\end{table}

\paragraph{Mesh Graph}

The problem of partitioning a  graph arises in
the context of computing a good ordering for the parallel
factorization of sparse, symmetric matrices \cite{psl90}. For load
balancing it is important to partition the graph in subgraphs of
approximately equal size. We constructed a small square grid of size
$32 \x 32$, and  partitioned the graph into 6 subgraphs.
In order to make the first 7 eigenvalues distinct, we let unnormalized
adjacency matrix be 
\[
  B=B_0 \otimes I + 0.7 (I \otimes B_0), 
\]
where $B_0$ is tridiagonal with 0 on the diagonal and 1 on the sub-
and superdiagonals. 
The results
are presented in  Table \ref{tab:mesh}.
\begin{table}[htb]
  \centering
  \begin{tabular}{lcccc}
    &     CPQR  & SSO &  QR+SSO &  Kmeans \\
    \hline
    Residual  &   0.364&  0.752&  0.361\\
    $\| \nabla_Q\|$ &  0.106 & $0.02\cdot 10^{-5}$&  $0.3\cdot 10^{-5}$\\
    Iterations&&   4/2 & 1/3\\
    NCut      &   0.270 &  0.585 &  0.270 &  0.350\\
    $\PsiCut$ & 0.131&  0.296&  0.131 &  0.180\\
    $\PhiCut$ &   0.0539&  0.262&  0.0539\\
                                 \hline
      \end{tabular}
      \caption{Mesh graph. First 7 eigenvalues:
        $1, 0.998,  0.997,  0.995,$ $0.992, 0.989,  0.988$.
    $L_k=0.0155$.  
\label{tab:mesh}  }
\end{table}

\begin{figure}[htbp!]   
	\begin{tabular}{ccc}
		\includegraphics[width=4cm, height=4cm]{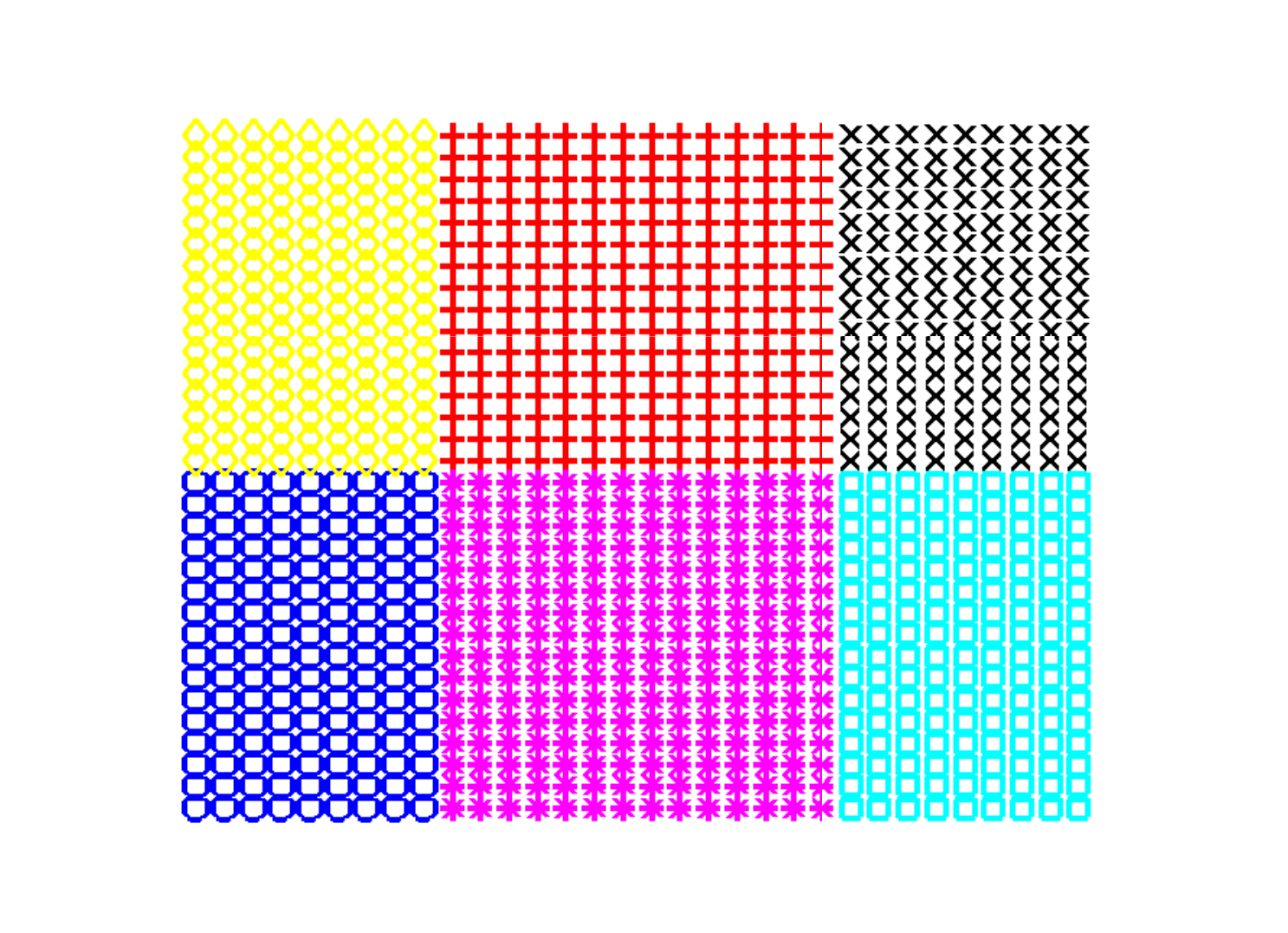} &
                                                                      \includegraphics[width=4cm, height=4cm]{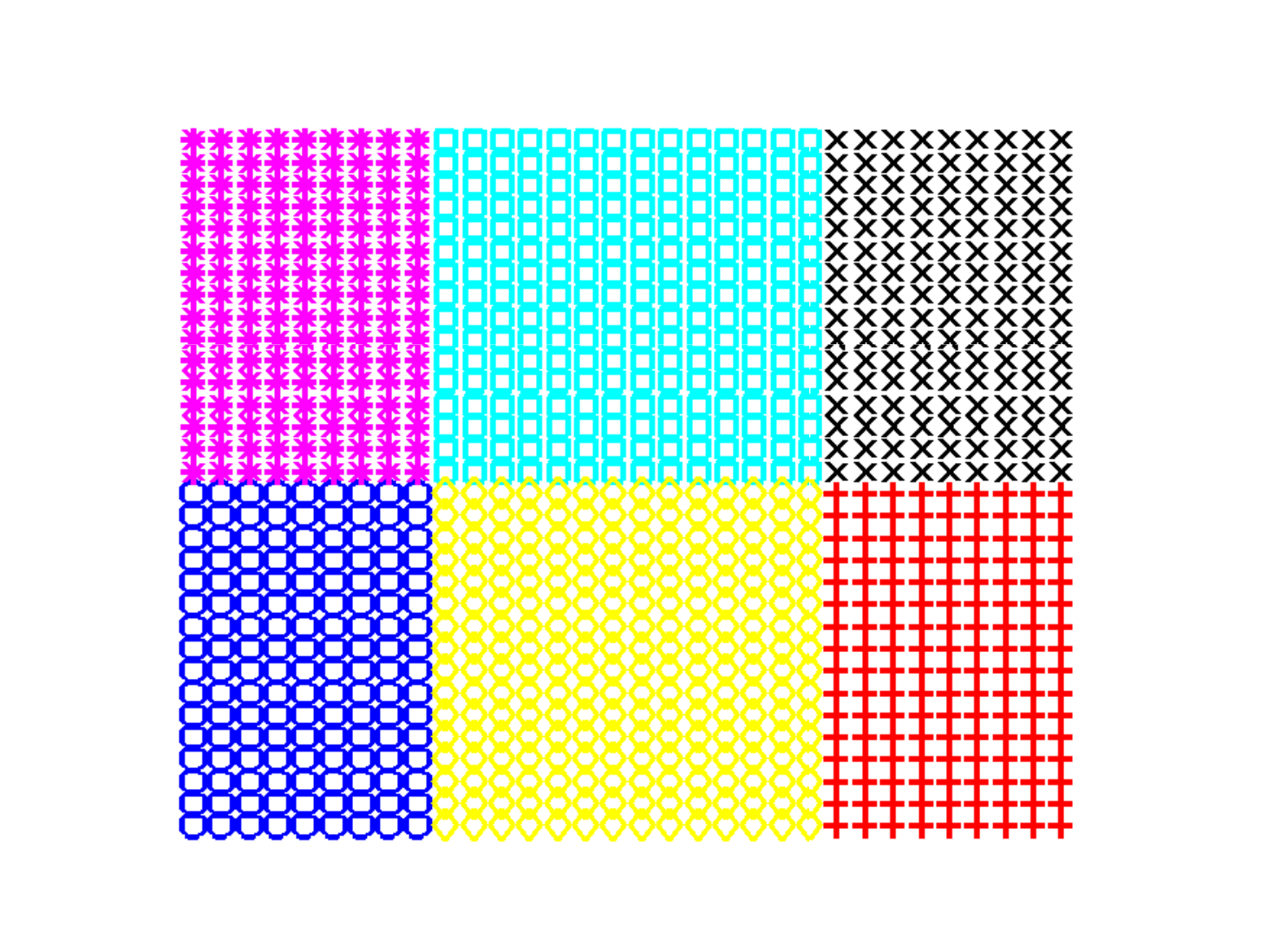}&\includegraphics[width=4cm, height=4cm]{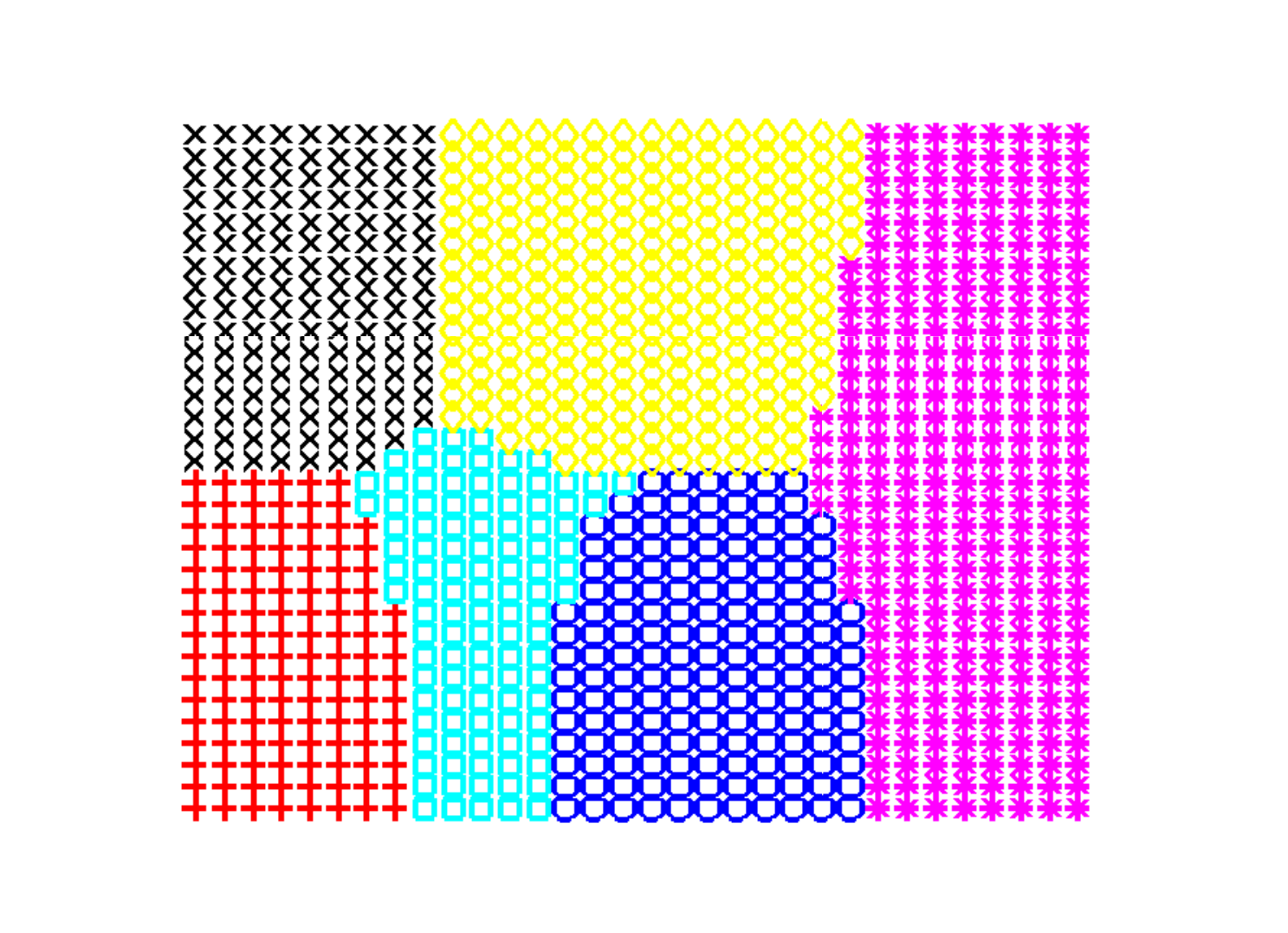} 
	\end{tabular}
	\caption{Mesh graph. Left to right: regions
          computed by CPQR, QR+SSO, and Kmeans. \label{fig:mesh}}  
\end{figure}

 SSO did not 
converge to the same point as QR+SSO. 
CPQR and QR+SSO gave the same partitionings, see
Figure \ref{fig:mesh}. The behavior of Kmeans was 
inconsistent; often it   gave
the same partitioning as  CPQR and QR+SSO.

\paragraph{Discussion of Examples}

We emphasize that most of the examples are chosen so that $L_k$ is small or at
least not large. This means that there is a strong  motivation for using
the spectral approach.

The first conclusion that can be drawn from the experiments is that
Kmeans is outperformed by the  other methods, because of its
inconsistent behavior and the fact that it never gives solutions of
better quality.

SSO never performs  better than
QR+SSO. CPQR gives a very good starting approximation for QR+SSO. For
most examples the difference between results of CPQR and QR+SSO is quite
small.   Since
the execution times of both methods are very low, it is  worth
the small effort in performing the SSO iterations, to get a
smaller residual, a smaller value of $\| \nabla_Q \|$,  and, in some
cases,  a partitioning of higher quality.

The algorithm QR$+$SSO  is fast: for  our largest example the
execution time was considerably shorter than the time for solving the
eigenvalue problem (using the Matlab function \texttt{eigs}). 
For our examples the number of iterations for QR+SSO is surprisingly small,
in spite of the fact that alternating iterations have linear
convergence rate. Therefore, for this type of problems,  it does not
seem worth the effort to 
develop a more advanced method for solving $\min \| X - \Psi Q^T\|$.  

 The computed
indicator matrix $\Phi$ was nonnegative for almost all test examples: 
for the yeast data  8 out of 1458 non-zeros were small and
negative.

\section{ Conclusions}
\label{sec:conclusions}

As we noted in the introduction the standard approach in spectral
graph partitioning is to start with a particular cut function and then
show that a relaxation leads to the solution of an eigenvalue problem
for the graph Laplacian. In this paper we do the converse: we start
with the eigenvalue/eigenvector problem, and show that if we restrict
the admissible solution to a set of vectors in indicator form, then we
get cut functions. From a theoretical point of view,   the cut
functions establish the relation between the eigenvalue  and graph
partitioning problems. From a practical point of view, the cut functions
can be used
to measure the quality of the computed partitioning.

We give a simple, efficient and robust algorithm for computing
indicator vectors from eigenvectors. Our experiments show that for
problems that are not far being $k$-partitionable the new algorithm is
to be preferred over Kmeans.
It is of interest to investigate the applicability of the
algorithm to a wider range of problems, but that is beyond the
scope of this paper.

\bibliographystyle{plain}

\bibliography{/media/lars/ExtHard/WORK/forskning/BIBLIOGRAPHIES/general,/media/lars/ExtHard/WORK/forskning/BIBLIOGRAPHIES/LE-papers}

\appendix
\appendixpage


\section{Closeness of a Graph to being $k$-Partitionable}
\label{app:Lk}

Fiedler \cite{fied73} defined the second smallest eigenvalue of the
unnormalized Laplacian to be the algebraic connectivity of the
graph. This measure has the property that if an edge is added, then
the connectivity  cannot become smaller. 

\begin{figure}[htbp!]   
	\begin{tabular}{ccc}
		\includegraphics[width=3.8cm, height=3.8cm]{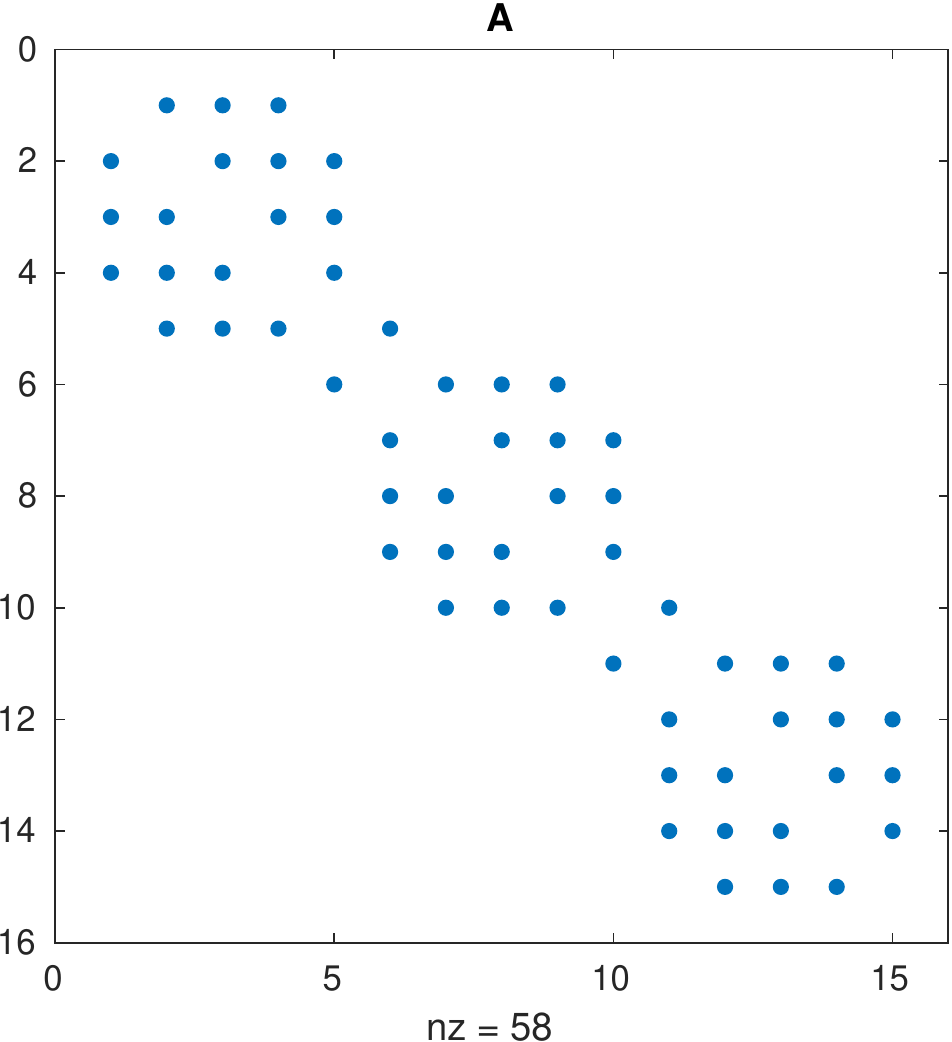} &
                                                                      \includegraphics[width=3.8cm, height=3.8cm]{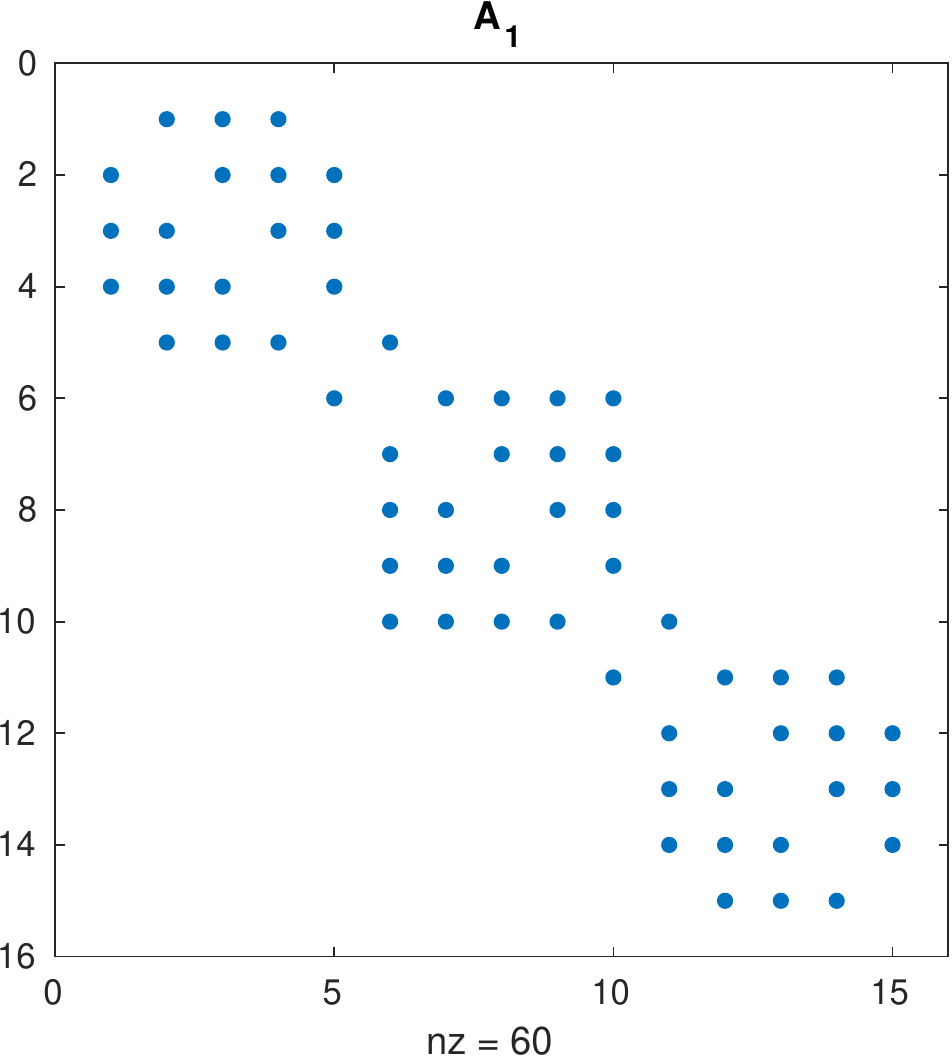}&\includegraphics[width=3.8cm, height=3.8cm]{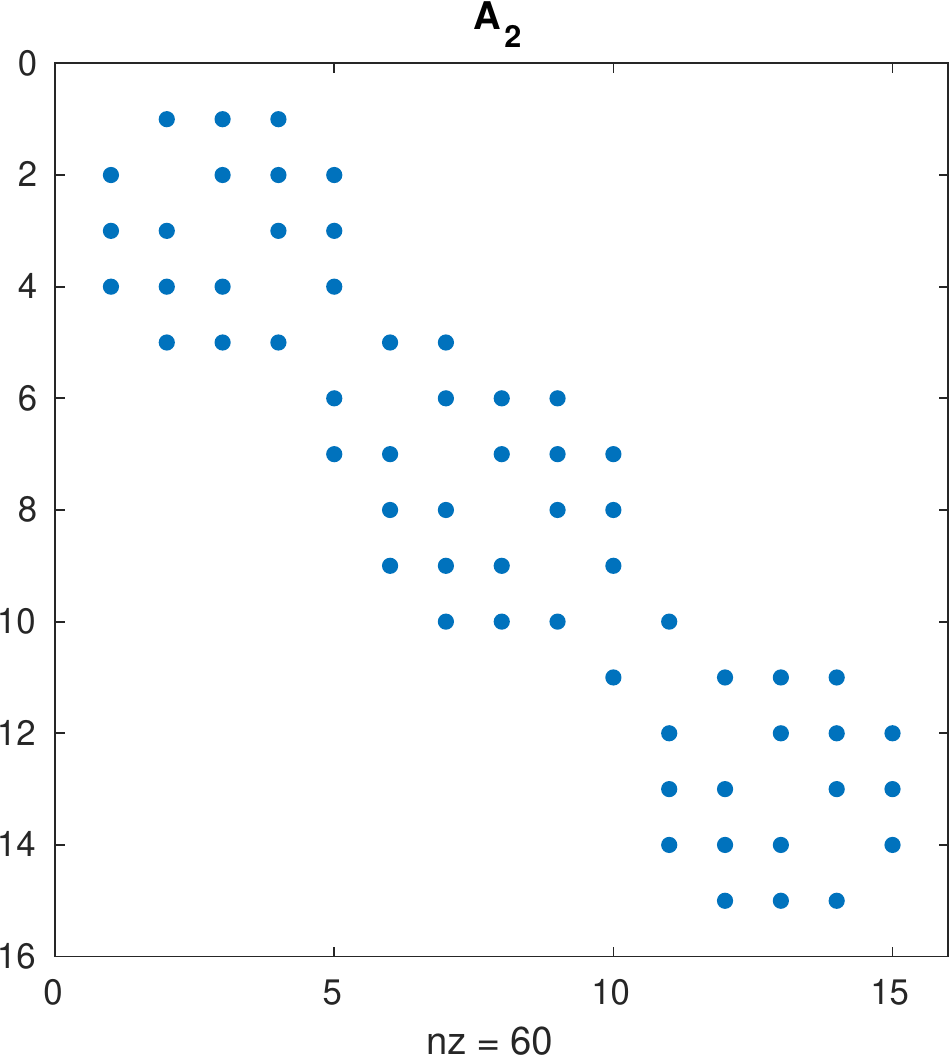} 
	\end{tabular}
	\caption{Adjacency matrices of three graphs  with
unweighted edges.\label{fig:three}}  
\end{figure}

Our measure $L_k$ for $k$-partitioning does not have the same
monotonicity property: if an edge is added, then usually $L_k$  becomes
larger, but it can also become smaller, 
depending on the structure. To illustrate this, we computed the
smallest eigenvalues of the normalized Laplacian of three graphs,
constructed so that they can naturally be partitioned into three
subgraphs, see Figure \ref{fig:three}.

%
%
%
In $A_1$ we added an edge centrally in the subgraph, far from the
partitioning positions, whereas in $A_2$ we added an edge close to one
of the partitioning positions. In Table \ref{tab:Lk-eig} we give the
smallest eigenvalues, cut functions, and $L_3$ for the three graphs.
\begin{table}[htb]
  \centering
  \begin{tabular}{lccc}
    &     $A$   & $A_1$  &  $A_2$ \\
    \hline
$1-\lambda_4$&0.802&    0.899  &  0.798\\
$1-\lambda_3$& 0.111&   0.106 &   0.152\\
$1-\lambda_2$  & 0.0335 &   0.0367&    0.0414\\
    $1-\lambda_1$&   0 & 0 & 0\\
    \hline
    $\PsiCut$ & 0.129 & 0.123 & 0.192\\
    $\PhiCut$ & 0.133 & 0.126 & 0.189\\
    $L_3$  & 0.116 &    0.112&    0.157\\
%
 \hline   
      \end{tabular}
  \caption{Three graphs. Four smallest eigenvalues $1-\lambda_i$ of the normalized
    Laplacian, cut functions,  and $L_3$.
    \label{tab:Lk-eig}  }
\end{table}
%
%
%
In $A_1$ we made the second subgraph slightly ``heavier'', which makes
it ``cheaper'' to $3$-partition the graph (cf. the cut functions). In
$A_2$ another edge will have to be broken in the partitioning into
subgraphs, and
therefore the distance of the third graph to being $3$-partitionable
is larger.

\section{Gradients}
\label{app:gradients}

\paragraph{Proof of Proposition \ref{prop:gradients}}

Assume that the indicator structure of $\Psi$ is fixed, and  write 
\[
  r(\Psi,Q) = \frac{1}{2} \| X- \Psi Q^T \|^2  =
  \frac{1}{2} \langle Y - \Psi, Y- \Psi \rangle =
    \frac{1}{2} \langle Y - \Psi, R \rangle,
\]
where $Y=XQ$, and $R=Y-\Psi$. Put $\Psi(t) =
\Psi + t \delta \Psi$, where $\delta \Psi$   has the same indicator
structure as $\Psi$. Then
\[
  \left.\frac{dr}{dt}\right|_{t=0} = - \langle \delta \Psi, R \rangle
  = \tr(R^T \delta\Psi)= 0,
\]
since, from \eqref{eq:Psicomp-3}-\eqref{eq:Psicomp-1}, the
corresponding columns of $\delta \Psi$ and $R$ have complementary
indicator structure.

To compute the gradient with respect to $Q$, we must take into account
that $Q$ is orthogonal. We parameterize $Q(t)$
along a geodesic curve in $\cO(k)$, starting from $Q(0)$ in the
direction $T$, where $T$ is 
a skew-symmetric matrix. It is well-known, see e.g. \cite{eas98},  that 
the  geodesic in $\cO(k)$ and its  $t$-derivative are
\[ 
  Q(t) = Q(0) \exp(tT), \qquad \dot{Q}(t)=Q(0)T \exp(tT).
\]
 It is no restriction to assume that $Q(0)=I$. 
 With $r(t) = \frac{1}{2} \langle Y - \Psi Q(t)^T, R \rangle$,
we have
\[
r'(0) = \left. \frac{dr}{dt}\right|_{t=0} = \langle  \Psi T, R \rangle =
  \tr\left(T^T \Psi^T R\right).
\]
 Let
\[
  T = S - S^T,
\]
where $S$ is strictly lower triangular, put $C=\Psi^TR $ and
consider $r'(0)= \tr(S^T C - SC)$.

We now use induction over the dimension $k$, and denote 
$S_k \in \RR^{k \x k}$ and $C_k \in \RR^{k \x k}$.
For $k=2$ we have
\begin{align*}
  \tr(S_2^T C_2 - S_2 C_2) &=
  \tr\left(
    \begin{pmatrix}
      0 & s_{21}\\
      0 & 0
    \end{pmatrix}
    \begin{pmatrix}
      c_{11} & c_{12}\\
      c_{21} & c_{22}
    \end{pmatrix}
    -
    \begin{pmatrix}
      0 & 0\\
      s_{21} & 0
    \end{pmatrix}
    \begin{pmatrix}
      c_{11} & c_{12}\\
      c_{21} & c_{22}
    \end{pmatrix}\right)\\
  &= (c_{21}-c_{12})s_{21}.
\end{align*}
Partition
\[
  S_{k} =
  \begin{pmatrix}
    S_{k-1} & 0\\
    s_{k}^T& 0
  \end{pmatrix}, \qquad
  C_{k}=
  \begin{pmatrix}
    C_{k-1} & c_{k}\\
    \bar c_{k}^T & *
  \end{pmatrix};
\]
(irrelevant elements are denoted by $*$). Straightforward computation gives,
\[
  S_{k}^T C_{k} - S_{k} C_{k} =
  \begin{pmatrix}
    S_{k-1}^T C_{k-1} + s_{k} \bar c_{k}^T & *\\
    * & 0 
  \end{pmatrix}
- \begin{pmatrix}
  S_{k-1} C_{k-1} & *\\
  * & s_{k}^T c_{k}
\end{pmatrix}.
\]
Thus
\begin{equation}\label{eq:tr-k+1}
  \tr\left(  S_{k}^T C_{k} - S_{k} C_{k} \right)=
  \tr\left(  S_{k-1}^T C_{k-1} - S_{k-1} C_{k-1} \right) +
                                 (\bar c_{k}^T - c_{k}^T)
                                 s_{k}. 
\end{equation}
By induction
\[
  \tr\left(  S_{k}^T C_{k} - S_{k} C_{k} \right)=
  \sum_{i=2}^k (\bar c_i^T - c_i^T) s_i,
\]
which shows that we can identify the operator $\nabla_Q$ with the
matrix
\[
  \cC =
  \begin{pmatrix}
    (\bar c_2^T - c_2^T) & (\bar c_3^T - c_3^T) & \cdots & (\bar c_k^T - c_k^T) 
  \end{pmatrix}.
\]
It is easy to see that
\[
  \| \nabla_Q \| = \| \cC \| = \frac{1}{2} \| C - C^T \|.
\]

\end{document}